\documentclass{article}
\usepackage{amsmath,amssymb}

\def\qed{\hspace*{\fill}{\Large$\square$}}
\def\bproof{$[\kern-.3ex[\,$}
\def\eproof{$\,]\kern-.3ex]$}

\newtheorem{lemma}{Lemma}
\newtheorem{theorem}{Theorem}

\title{A Characterization of Similarity Maps Between Euclidean Spaces Related
  to the Beckman--Quarles Theorem}

\date{\normalsize Abbreviated title: Similarity Maps Between Euclidean Spaces}

\author{Jobst Heitzig\\[1mm]
  \small Institut f\"ur Mathematik, Universit\"at Hannover\\ 
  \small Welfengarten 1, 30167 Hannover, Germany\\
  \small {\tt heitzig@math.uni-hannover.de}}

\begin{document}
\maketitle

\begin{abstract}
  \noindent
  It is shown that each continuous transformation $h$ from Euclidean $m$-space
  ($m>1$) into Euclidean $n$-space that preserves the equality of distances
  (that is, fulfils the implication
  $|x-y|=|z-w|\Longrightarrow|h(x)-h(y)|=|h(z)-h(w)|$) is a similarity map. 
  The case of equal dimensions already follows from the Beckman--Quarles
  Theorem.

  ~


  ~


\end{abstract}

\section*{Introduction}
\noindent
In the theory of metric spaces or, more generally, distance sets and
distance spaces (cf.\ \cite{u-Hei01,u-Hei02}), 
distance equality preserving transformations, that is, those that
preserve equations of the form $d(x,y)=d(z,w)$, are far more general 
than transformations which preserve all distances accurately ({\em isometries}) or
up to a multiplicative constant ({\em similarity maps}). 
Using monoid homomorphisms, it shall be shown here that in case of Euclidean
spaces of dimension at least two, both classes of maps are in fact the same
when continuity is required in addition. 

Throughout, $G_n$ will denote the
group of motions (=self-isometries) of Euclidean $n$-space $\mathbb R^n$,
equipped with the topology of pointwise convergence. Double square brackets
\bproof\dots$\!$\eproof\ are used to enclose proof details.

\section*{A simple fact about distance functions}

Although I will apply it only for the case of Euclidean metric, the
following lemma states a simple correspondence between distance equality
preserving maps and monoid homomorphisms in presence of a very general kind of
distance function. Formally, a distance function is here just a map $d:X^2\to
M$ into some set $M$, and a transformation $h:(X,d)\to\mathbb R^n$ is {\em
distance equality preserving} if and only if the implication
$d(x,y)=d(z,w)\Longrightarrow|h(x)-h(y)|=|h(z)-h(w)|$ holds for all
$x,y,z,w\in X$. 

\begin{lemma}\label{lem:eqpreshom}
Let $(X,\cdot,1)$ be a monoid, 
$n\geqslant 0$, and $d$ a left translation-invariant
distance function on $X$, that is, with $d(zx,zy)=d(x,y)$. 
If $h:(X,d)\to\mathbb R^n$ is a distance equality preserving map, there is
a monoid homomorphism $f:(X,\cdot,1)\to G_n$ such that
$h(x)=f(x)\big(h(1)\big)$ for all $x\in X$.

If, moreover, $h[X]$ is not contained in any affine hyperplane of\/ $\mathbb
R^n$, this $f$ is unique.
\end{lemma}
{\em Proof.\/}
Let us first consider the ``non-degenerate'' case where $h[X]$ is not
contained in a hyperplane. Then there are $n+1$ points 
$x_0,\dots,x_n\in X$ such that the set $\{h(x_i):i\in n+1\}$ is not contained 
in any hyperplane. For each $x\in X$, there is a unique motion
$f_x$ of $\mathbb R^n$ with $f_xh(x_i)=h(xx_i)$ for all $i$ 
\bproof Since $d$ is left translation invariant and $h$ is distance
equality preserving, $|h(xx_i)-h(xx_j)|=|h(x_i)-h(x_j)|$, hence also 
$\{h(xx_i) | i\in n+1\}$ is not contained in a hyperplane. 
Thus, for all $a\in\mathbb R^n$, there is a unique $b\in\mathbb R^n$ such that
$|b-h(xx_i)|=|a-h(x_i)|$ for all $i$. Put $f_x(a):=b$. 
In particular, $f_x(h(x_i))=h(xx_i)$ for
all $i$. Since for all $a,a'\in\mathbb R^n$, the distance $|a-a'|$ 
is a function of the two $(n+1)$-tuples 
of distances $\big(|a-h(x_i)|\big)_i$ and $\big(|a'-h(x_i)|\big)_i$, 
$f_x$ is a motion\eproof. 
Now $h(xy)=f_x(h(y))$ for all $x\in X$ 
\bproof $f_x(h(y))$ is the unique $b\in Y$ with
$|b-h(xx_i)|=|h(y)-h(x_i)|$ for all $i$,  
and $h(xy)$ is such a $b$\eproof, 
in particular $h(x)=f_x(h(1))$, and the map $f:x\mapsto f_x$ is a monoid
homomorphism 
\bproof $f_x\circ f_y$ is a motion with $f_x\circ
f_y(h(x_i))=f_x(h(yx_i))=h(xyx_i)$ for all $i$, hence it equals 
$f_{xy}$\eproof.

On the other hand, let $g:x\mapsto g_x$ be another homomorphism with
$h(x)=g_x(h(1))$ for all $x$. Then $g_x(h(x_i))=g_x\circ
g_{x_i}(h(1))=g_{xx_i}(h(1))=h(xx_i)$ 
for all $x$ and $i$, hence $g_x=f_x$ for all $x$, that is, $g=f$.

For the degenerate case, let $h':=i^{-1}\circ h$, 
where $i$ is an isometry between some $\mathbb R^k$ (with $k<n$) 
and the affine hull of $h[X]$. Then $h':(X,d)\to\mathbb R^k$ 
is distance equality preserving and
``non-degenerate'', hence there is a corresponding monoid homomorphism
$f':(X,\cdot,1)\to G_k$. 
Moreover, there is an embedding $g:G_k\to G_n$ 
such that $g(m)\circ i=i\circ m$ for all motions $m$ of $\mathbb R^k$.
Then $f:=g\circ f'$ is a homomorphism such that 
\[ f(x)\big(h(1)\big)=\big(g(f'_x)\circ i\big)\big(h'(1)\big)=
   (i\circ f'_x)\big(h'(1)\big)=i\circ h'(x)=h(x). \]
\qed

~

Moreover, it is easily verified that in case of $X=\mathbb R^m$, either none
or both of $h$ and $f$ are continuous.

\section*{The case of a one-dimensional domain}

It is well known that each continuous representation of the (additive) group 
$\mathbb R$ by motions of $\mathbb R^n$ is of the following form.
Let $k$ be a motion of $\mathbb R^n$ 
and $s$ a non-negative integer with $2s\leqslant n$.
For each $i\in\{1,\dots,s\}$, let $E_i$ be the plane
spanned by the standard unit vectors $e_{2i-1}$ and $e_{2i}$, and 
$\alpha_i>0$.
Moreover, let the vector $b\in\mathbb R^n$ be orthogonal to all $E_i$. 
Then the map $f:x\mapsto f_x$ defined by 
\[\textstyle
  A_i:=\left(\textstyle\addtolength{\arraycolsep}{-.5ex}\begin{array}{rr}
                \scriptstyle \cos(x\alpha_i) & \scriptstyle -\sin(x\alpha_i) \\
                \scriptstyle \sin(x\alpha_i) & \scriptstyle \cos(x\alpha_i)\end{array}\right)\text{~~and~~}
  k\big(f_x(a)\big):=xb+\left(\textstyle\addtolength{\arraycolsep}{-.5ex}\begin{array}{cccc}
        \framebox{$\scriptstyle A_1$}\\
        &\scriptstyle\ddots \\
        &&\framebox{$\scriptstyle A_s$}\\
        &&&\framebox{$\scriptstyle 1_{n-2s}$}
        \end{array}\right)\cdot k(a) \]
is a continuous group homomorphism from $\mathbb R$ into $G_n$, where
$1_{n-2s}$ is a unit matrix of $n-2s$ dimensions.
Each $f_x$ is a composition of rotations in the planes $k^{-1}[E_i]$ with 
centre $k^{-1}(0)$ and angles $x\alpha_i$, and of a translation perpendicular
to all those planes. 

Obviously, the continuous distance equality preserving maps from $\mathbb R^1$
to $\mathbb R^n$ are exactly the ``generalized helices'' $h(x)=f_x(a)$ with
$f$ of the form just described and arbitrary $a\in\mathbb R^n$.

\section*{The case of a higher-dimensional domain}
We will now see that between Euclidean spaces of higher dimension than one,
the continuous distance equality preserving maps are already similarity maps.

\begin{lemma}\label{lem:uniqueangles}
  For a function $C:[0,\infty)\to[0,\infty)$ of the form 
\[C(r)=(r\lambda)^2+\sum_{k=1}^s
  \left(\lambda_k\sqrt{2-2\cos(r\kappa_k)}\right)^2 \]
  with $0<\kappa_1<\cdots<\kappa_s$ and $\lambda_k>0$, all coefficients
  $\kappa_k,\lambda_k$
  are uniquely determined.
\end{lemma}
{\em Proof.}
For $n>0$, the $(4n+1)$-st derivative of $C$ fulfils
\[ \frac{C^{(4n+1)}(r)}{2\kappa_s^{4n+1}} = \lambda_s^2\sin(r\kappa_s) +
\sum_{k=1}^{s-1}\lambda_k^2\Big(\frac{\kappa_k}{\kappa_s}\Big)^{4n+1}\sin(r\kappa_k)
\quad\to\quad\lambda_s\sin(r\kappa_s)\] 
for $n\to\infty$.
Therefore, $\kappa_s$ is the smallest $\mu>0$ with $C^{(4n+1)}(r)=O(\mu^{4n+1})$
for almost all $r$, and $\lambda_s^2=\sup_r\lim
C^{(4n+1)}(r)/2\kappa_s^{4n+1}$. Subtracting the $s$-term, one can now
inductively determine all $\kappa_k$ and $\lambda_k$.
\qed

\begin{theorem}\label{th:sims}
The continuous distance equality preserving maps from\/ $\mathbb R^m$ to\/
$\mathbb R^n$ with $m>1$ are exactly the similarity maps.
\end{theorem}
{\em Proof.}
Since a map $h:\mathbb R^m\to\mathbb R^n$ ($m>1$) is a similarity if and only
if $h|_E$ is one for all affine planes $E\subseteq\mathbb R^m$,
we may assume that $m=2$. Let
$h:\mathbb R^2\to\mathbb R^n$ be continuous and distance equality preserving
with $h(0,0)=z$.
From Lemma \ref{lem:eqpreshom} we know that $h$ is of the form
\[ h(x,y)=f_{(x,y)}(z)=\varphi_x\big(\psi_y(z)\big)=\psi_y\big(\varphi_x(z)\big), \]
where $\varphi:x\mapsto\varphi_x$ and $\psi:y\mapsto\psi_y$
are continuous group
homomorphisms from $\mathbb R^2$ into $G_n$. The motions $\varphi_x$ and
$\psi_y$ are of the form
\[ \varphi_x(v)=a_x+A_xv\quad\text{and}\quad\psi_y(v)=b_y+B_yv, \]
where $A_x$ and $B_y$
are orthogonal matrices with $A_xB_y=B_yA_x$
\bproof since
$a_x+A_xb_y+A_xB_yv=\varphi_x\big(\psi_y(v)\big)=\psi_y\big(\varphi_x(v)\big)=b_y+B_ya_x+B_yA_xv$
for all $v\in\mathbb R^n$ implies $a_x+A_xb_y=b_y+B_ya_x$, 
hence $A_xB_yv=B_yA_xv$ for all $v$, that is, $A_xB_y=B_yA_x$\eproof.

Therefore, there is a unitary complex matrix $P$ 
such that $A'_x:=P^{-1}A_xP$ and $B'_y:=P^{-1}B_yP$
are (complex) diagonal matrices for all $x,y$
\bproof Choose $\xi,\eta\in\mathbb R$ so that $A_\xi$ and $B_\eta$ 
have a minimal number of real eigenvalues among all $A_x$ resp.\ $B_y$. 
The commuting orthogonal matrices $A_\xi$ and $B_\eta$ have a common
diagonalization $P^{-1}A_\xi P$, $P^{-1}B_\eta P$ with some unitary $P$.
Since for all $x,y\in\mathbb R$, the planes of rotation of $A_x$ and $B_y$ are among 
those of $A_\xi$ and $B_\eta$, respectively, every complex eigenvector of
$A_\xi$ or $B_\eta$ is an eigenvector of $A_x$ or $B_y$, respectively.
Hence also $A'_x$ and $B'_y$ are diagonal\eproof.
Now $x\mapsto A'_x$ and $y\mapsto B'_y$ are continuous homomorphisms into the
group of unitary diagonal matrices, hence there are coefficients
$\alpha_j,\beta_j\in\mathbb R$ such that 
\[ A'_x=\mbox{\rm diag}(e^{ix\alpha_1},\dots,e^{ix\alpha_n})\text{~~and~~}
   B'_y=\mbox{\rm diag}(e^{iy\beta_1},\dots,e^{iy\beta_n}). \]

The map $c:[0,\infty)\to[0,\infty)$, 
$|x-y|\mapsto |h(x)-h(y)|$ is well-defined since $h$ is distance equality
preserving, and fulfils  
\[ c(r)=|h(0,0)-h(r\cos\gamma,r\sin\gamma)|
       =\big|z-\varphi_{r\cos\gamma}\big(\psi_{r\sin\gamma}(z)\big)\big| \]
for all $r\geqslant 0$ and $\gamma\in\mathbb R$.
Note that $\varphi_{r\cos\gamma}\circ\psi_{r\sin\gamma}$ is a motion whose
matrix $A_{r\cos\gamma}B_{r\sin\gamma}$ has the complex eigenvalues
$e^{ir((\cos\gamma)\alpha_j+(\sin\gamma)\beta_j)}$, 
$j=1\dots n$. Using elementary geometry, we see that
\[ c(r)^2=\big(r\lambda(\gamma)\big)^2+\sum_{k=1}^{s(\gamma)}
   \left(\lambda_k(\gamma)\sqrt{2-2\cos\big(r\kappa_k(\gamma)\big)}\right)^2
   \quad\text{for all~}\gamma\in\mathbb R, \] 
with $\lambda_k(\gamma)>0$ and $\kappa_k(\gamma)>0$ for all $k$. 
Indeed, $r\lambda(\gamma)$ is the length
of the translational part of $\varphi_{r\cos\gamma}\circ\psi_{r\sin\gamma}$;
each $\lambda_k(\gamma)$ is a radius of rotation for some rotational part of
$\varphi_{r\cos\gamma}\circ\psi_{r\sin\gamma}$, that is, the distance from $z$
to the 
affine $(n-2)$-dimensional subspace of $\mathbb R^n$ which is fixed under that
rotation; and $r\kappa_k(\gamma)$ is the corresponding angle of that rotation.

Note that, by definition of the $\alpha_j,\beta_j$, each $\kappa_k(\gamma)$ is
of the form $|(\cos\gamma)\alpha_j+(\sin\gamma)\beta_j|$ for some $j$.
But because of Lemma \ref{lem:uniqueangles}, the set
$K:=\{\kappa_k(\gamma) : k=1\dots s(\gamma)\}$ can be determined from $c$
and is thus the same for all $\gamma\in\mathbb R$. Assume that $\kappa\in K$. Then
for each $\gamma\in\mathbb R$, there is $j\in\{1,\dots,n\}$ with
$|(\cos\gamma)\alpha_j+(\sin\gamma)\beta_j|=\kappa$. This is only possible if
$\kappa=\alpha_j=\beta_j=0$ for some $j$, in contradiction to
$\kappa>0$. Hence $K$ is empty and $c$ is linear, which means that $h$ is a
homometry.\qed

~

The special case of equal dimensions $n=m$ also follows from the 
Beckman--Quarles Theorem \cite{u-BQ53} 
which says that a map $f:\mathbb R^n\to\mathbb R^n$ ($n\geqslant
2$) with $e(x,y)=1\Longrightarrow ef(x,y)=1$ is already an isometry. 
It is not obvious, however, if one of its proof 
(see, e.\thinspace g.~\cite{u-Len91,u-Les95}) or variants (cf.~\cite{u-Ben80,u-Ben82,u-Hot97,u-Huc86,u-Rad84,u-Vas86})
generalizes to the case of arbitrary different dimensions.

{
\small
\def\NDASH{--}
\bibliographystyle{amsplain} 
\bibliography{../all.bib}
}

\end{document}